\documentclass[12pt]{article}
\usepackage{amsmath,amssymb}
\usepackage{mathrsfs}
\usepackage{mathtools}
\usepackage{amsthm}
\usepackage{amscd}
\usepackage[all]{xy}
\usepackage{enumitem}
\usepackage[top=25mm,bottom=25mm,left=20mm,right=20mm]{geometry}
\usepackage{comment}

\newcommand{\rmF}{\mathrm{F}}

\newcommand{\rmH}{\mathrm{H}}
\newcommand{\frakm}{\mathfrak{m}}

\newcommand{\calO}{\mathcal{O}}

\newcommand{\frakp}{\mathfrak{p}}
\newcommand{\Q}{\mathbf{Q}}
\newcommand{\frakq}{\mathfrak{q}}

\newcommand{\Z}{\mathbf{Z}}
\theoremstyle{plain}
\newtheorem{definition}{Definition}
\newtheorem{thm}[definition]{Theorem}
\newtheorem{prop}[definition]{Proposition}

\newtheorem{lemma}[definition]{Lemma}
\newtheorem{remark}[definition]{Remark}

\newtheorem{cor}[definition]{Corollary}

\numberwithin{definition}{section}
\numberwithin{equation}{section}

\DeclareMathOperator{\Dt}{Dt}

\DeclareMathOperator{\Hom}{Hom}

\DeclareMathOperator{\Spec}{Spec}

\DeclareMathOperator{\Sw}{Sw}
\DeclareMathOperator{\ab}{ab}
\DeclareMathOperator{\ch}{char}

\DeclareMathOperator{\dt}{dt}
\DeclareMathOperator{\fil}{fil}

\DeclareMathOperator{\ord}{ord}

\DeclareMathOperator{\rsw}{rsw}
\DeclareMathOperator{\sw}{sw}

\newcommand{\map}[3]{#1\colon#2\to#3}
\newcommand{\lra}{\longrightarrow}
\usepackage[%
 setpagesize=false,%
 bookmarks=true,%
 bookmarksdepth=tocdepth,%
 bookmarksnumbered=true,%
 colorlinks=false,%
 pdftitle={},%
 pdfsubject={},%
 pdfauthor={},%
 pdfkeywords={}%
]{hyperref}
\begin{document}
\title{The abelian fundamental group with modulus\\ in mixed characteristic}
\author{Ryosuke Ooe}
\date{}

\maketitle
\abstract{We define the abelian fundamental group with modulus of a regular flat scheme over a discrete valuation ring, taking into account wild ramification along a divisor. Our definition provides a mixed-characteristic analogue of the abelian fundamental group with modulus introduced by Kerz--Saito for smooth schemes over a perfect field. In this setting, we prove a Lefschetz-type theorem for strictly semi-stable schemes: restriction to a hypersurface of sufficiently large degree relative to the ramification induces an isomorphism of the abelian fundamental groups. } 
\section*{Introduction}

Let $k$ be a field, and let $X$ be a smooth projective connected scheme over $k$. 
Fix a projective embedding and let $Y$ be a smooth hypersurface. 
Grothendieck \cite{SGA2} proved that if $\dim X \ge 3$, then the immersion induces an isomorphism
\begin{equation}\label{0.1}
\pi_1(Y) \xrightarrow{\cong} \pi_1(X)
\end{equation}
of \'etale fundamental groups. 

If we consider quasi-projective schemes, one must take into account ramification along a divisor, since the map \eqref{0.1} is not an isomorphism in positive characteristic. 
In the tame case, this was settled by Esnault--Kindler \cite{EK16}, who proved a Lefschetz-type theorem for the tame fundamental group. In the wild case, for a smooth scheme $X$ over a perfect field $k$ of positive characteristic and a simple normal crossing divisor $D$ on $X$, 
Kerz--Saito \cite{KS14} defined the abelian fundamental group with modulus, denoted by $\pi_1^{\ab}(X,D)$. Let $U$ be the complement of $D$. It can be understood as the quotient of $\pi_1^{\ab}(U)$ classifying \'etale covers of $U$ with ramification bounded by $D$, while it is defined as the dual of the subgroup of characters of $\pi_1^{\ab}(U)$ whose ramification is bounded by $D$.

There are at least two approaches to bounding the wild ramification of $\ell$-adic sheaves along divisors. 
One approach uses classical ramification theory for local fields with perfect residue fields, 
via restriction of sheaves to curves. 
Alternatively, one can use the ramification theory of Abbes--Saito \cite{AS02}, 
which applies to local fields at the generic points of the irreducible components of divisors. 

In \cite{KS14}, Kerz--Saito established the following Lefschetz-type theorem for the abelian fundamental group with modulus. 

\begin{thm}[{\cite[Theorem 1.1]{KS14}}]
    Assume that $Y$ is sufficiently ample for $(X,D)$. Then 
\[
\pi_1^{\ab}(Y,E) \lra \pi_1^{\ab}(X,D) 
\]
is an isomorphism if $\dim X \ge 3$ and a surjection if $\dim X=2$. 
\end{thm}
Here, the condition that $Y$ is sufficiently ample for $(X,D)$ is formulated in terms of the vanishing of certain cohomology groups. It is in particular satisfied when the degree of $Y$ is sufficiently large.

It is natural to ask whether the definition and the Lefschetz theorem of Kerz--Saito extend to schemes over a discrete valuation ring. 
The purpose of this paper is to address this question in the mixed characteristic case.
Let $U$ be a quasi-projective smooth scheme over a complete discrete valuation field $K$ of characteristic zero. We assume that there exists a strictly semi-stable projective model $X$ over $\calO_K$ such that $D=X-U$ is a simple normal crossing divisor.  

We define the abelian fundamental group $\pi_1^{\ab}(X, D)$ in the same way as Kerz--Saito. As in the equal characteristic case, we can bound wild ramification in two approaches. We prove that the two notions of bounded ramification for characters are equivalent in Theorem \ref{1.2}. 

 Let $Y$ be a regular hypersurface flat over $\calO_K$ such that $Y$ satisfies suitable transversality conditions. We establish the following Lefschetz-type theorem. 
\begin{thm}[Theorem \ref{main}]
    Assume that $Y$ is sufficiently ample for $(X,D)$. Then 
\[
\pi_1^{\ab}(Y,E) \lra \pi_1^{\ab}(X,D) 
\]
is an isomorphism if $\dim X_K \geq 3$ and a surjection if $\dim X_K = 2$. 
\end{thm}

In order to formulate the condition that $Y$ is sufficiently ample for $(X,D)$, we make use of the module of Frobenius–Witt differentials, recently introduced by Saito \cite{Sa22}. This condition plays a role when deducing the theorem from its logarithimic version. The key fact for the proof of the logarithimic version is that the Swan conductor for characters has an upper bound (Proposition \ref{1.3}). This fact is true only in mixed characteristic. 

The paper is organized as follows. 
In Section 1.1, we begin by recalling the necessary background on ramification theory and review properties of the Swan conductor. 
In Section 1.2, we define the abelian fundamental group with modulus in two different ways. We prove that these definitions are equivalent. 
In Section 2, we define the condition that $Y$ is sufficiently ample for $(X, D)$. We prove the main theorem, establishing the Lefschetz-type result. 
\vskip\baselineskip

\textbf{Acknowledgment}

The author would like to express his sincere gratitude to his advisor, Professor Takeshi Saito, for carefully reading the drafts, pointing out errors, and helping with the proof of Proposition \ref{1.3}. The author is also grateful to Professor Shuji Saito for kindly answering his questions. 
This work was supported by JSPS KAKENHI Grant Number JP24KJ0833.

\tableofcontents
\setcounter{section}{0}
\section{Ramification theory}
\subsection{Local theory}
In this subsection, let $K$ be a henselian discrete valuation field with residue field $F$ of characteristic $p>0$. Let $\calO_K$ denote its valuation ring and $\frakm_K$ be the maximal ideal of $\calO_K$. Let $\pi$ be a uniformizer of $K$. Let $K_1$ be the fraction field of the henselization of the local ring of $\calO_K[T]$ at the prime ideal $\frakm_K\calO_K[T]$. 
\begin{definition}[{\cite[Definition (2.1)]{Ka89}}]
  Let $\chi\in \rmH^1(K,\Q/\Z)$ be a character. The Swan conductor of $\chi$ is defined as the smallest integer $n$ satisfying
  \[
  \{\chi, 1+\frakm_{K_1}^{n+1}\}=0
  \]
  in $\mathrm{Br} K_1$. 
\end{definition}

There is another definition of the Swan conductor using the upper ramification group defined by Abbes--Saito \cite{AS02}. The equality of two definitions of the Swan conductor was proved in \cite{KS19}. 

The total dimension $\dt(\chi)$ is defined as a non-logarithmic variant of the Swan conductor. We have either $\dt(\chi)=\sw(\chi)+1$ or $\dt(\chi)=\sw(\chi)$. For details, see \cite{Ma97}, \cite{Sa23}.

We say that $\chi$ is of type I if $\dt(\chi)=\sw(\chi)+1$ and of type II if $\dt(\chi)=\sw(\chi)$. If the residue field of $K$ is perfect, then $\chi$ is always of type I.

\begin{prop}\label{1.3}
  Assume that the characteristic of $K$ is $0$. Let $n$ be a positive integer. Then, there exists an integer $N_n$ such that we have $\sw \chi\le N_n$ for all characters $\chi$ of order $p^n$. In particular, we can take $N_1=\lceil ep/(p-1)\rceil$ where $e$ denotes the absolute ramification index of $K$. 
\end{prop}
\begin{proof}
  It suffices to show that there exists an integer $N_n$ such that we have 
  \[
  1+\frakm^{N_n+1}_{K_1}\subset (K_1^\times)^{p^n}
  \] 
 since we have $\{\chi, (K_1^\times)^{p^n}\}=0$ for characters $\chi$ of order $p^n$.
  
  For a positive integer $N$ and an element $a\in \frakm^{N+1}_{K_1}$, consider the following equation with coefficients in $\calO_{K_1}$:
  \begin{equation}\label{1.1}
  (1+\pi^iX)^{p^n}=1+a.
  \end{equation}
  Choose an integer $M_n$ such that, for every $2\le j\le p^n$, the inequality  
  \[
  ne+i<e\cdot \ord_p\binom{p^n}{j}+j\cdot i
  \]
  holds whenever $i>M_n$.
  By the henselian property, equation (\ref{1.1}) has a root if $i>M_n$ and $N\ge ne+i$. Hence, taking $N_n=ne+M_n$, we obtain
  \[
  1+\frakm_K^{N_n+1}\subset (K_1^\times)^{p^n}. 
  \]
  When $n=1$, we can take $M_1=\lceil e/(p-1)\rceil$, and thus $N_1=\lceil ep/(p-1)\rceil$. 
  
  \end{proof}

\begin{lemma}\label{1.4}
  Assume that the characteristic of $K$ is zero. Let $e$ be the absolute ramification index of $K$, and set $e'=ep/(p-1)$. Let $\chi\in H^1(K, \Q/\Z)$ be a character. Then we have $\sw(p\chi)\le \sw(\chi)$. Moreover, if $\sw(\chi)>e'$, we have $\sw(p\chi)=\sw(\chi)-e$.
\end{lemma}
\begin{proof}
The inequality $\sw(p\chi)\le \sw(\chi)$ follows directly from the definition of the Swan conductor. The latter assertion follows from \cite[Lemma 2.19]{BN23}. 
\end{proof}

As a refinement of the Swan conductor, Kato \cite{Ka89} defined the injection
\[
\rsw(\chi)\colon \frakm_K^{-n}/\frakm_K^{-n+1}\lra \Omega^1_{F}(\log)
\]
called the refined Swan conductor if $n=\sw(\chi)>0$. 

As a non-logarithmic variant of the refined Swan conductor, Saito \cite{Sa23} defined the injection
\[
\ch(\chi)\colon \frakm_K^{-m}/\frakm_K^{-m+1}\lra \rmH_1(L_{F^{1/p}/\calO_K})
\]
called the characteristic form if $m=\dt(\chi)>1$.
(The fact that the image of the characteristic form is contained in $\rmH_1(L_{F^{1/p}/\calO_K})$ is proved in \cite{O25}.) 

We also define $\ch'(\chi)$ to be the composition
\[
\ch'(\chi)\colon \frakm_K^{-m}/\frakm_K^{-m+1}\xrightarrow{\ch(\chi)}    \rmH_1(L_{F^{1/p}/\calO_K})\lra \rmH_1(L_{F^{1/p}/\calO_K})\otimes_{F^{1/p}}F
\]
where $F^{1/p}\to F$ denotes the $p$-th power map. 

\begin{lemma}\label{1.5}
 Assume that the characteristic of $K$ is zero. Let $\chi\in \rmH^1(K, \Q/\Z)$ be a character. If $\chi$ is of type \textup{II}, then the Swan conductor of $\chi$ is divisible by $p$. 
\end{lemma}
\begin{remark}
  This lemma is trivial when $\ch K=p$ by the definition of the total dimension due to Matsuda \textup{\cite{Ma97}}. 
\end{remark}
\begin{proof}
Let $n$ be the Swan conductor of $\chi$. By assumption, the total dimension of $\chi$ is also $n$. 

  Write
  \[
  \rsw(\chi)=(\alpha+\beta d\log \pi)/ \pi^n \in \frakm_K^{-n}/\frakm_K^{-n+1}\otimes \Omega^1_F(\log)
  \]
  where $\alpha\in \Omega^1_F, \beta \in F^{\times}$. 
By \cite[Lemma 2.17]{BN23}, we have $d\beta=-n\alpha$. 

  Similarly, write 
  \[
  \ch(\chi)=(\alpha'+\beta' w\pi)/ \pi^n \in \frakm_K^{-n}/\frakm_K^{-n+1}\otimes \rmH_1(L_{F^{1/p}/\calO_K})
  \]
  where $\alpha'\in \Omega^1_F, \beta' \in F^{\times}$. 

  Since we assume $\chi$ to be of type II, the refined Swan conductor is the image of the characteristic form under the map induced by $\rmH_1(L_{F^{1/p}/\calO_K})\to \Omega^1_F(\log)$ by \cite[Proposition 3.4]{O25}. Hence, we have $\beta=\pi\beta'=0$ and $-n\alpha =0$. The refined Swan conductor is non-zero, so $\alpha$ is non-zero. Hence, $n$ is divisible by $p$. 
\end{proof}

\begin{cor}\label{1.11}
Let $m$ be a positive integer prime to $p$. For $\chi\in \rmH^1(K,\Q/\Z)$, the following conditions are equivalent:

\textup{1}. $\sw(\chi)\le m-1$.

\textup{2}. $\dt(\chi)\le m$. 
\end{cor}
\begin{proof}
    The implication $1\Rightarrow 2$ is trivial. We show $2\Rightarrow 1$. When $\dt\chi\le m$, we have $\sw\chi\le m$. If we have $\sw(\chi)=m$, then $\dt(\chi)=m+1$ by Lemma \ref{1.5} since $m$ is prime to $p$. Hence, we have $\sw\chi\le m-1$. 
\end{proof}

\subsection{Global theory}
Let $X$ be a regular excellent scheme and let $C$ be a simple normal crossing divisor on $X$. Let $U$ be the complement of $C$. Let $\{C_i\}_{i=1,\dots r}$ be irreducible components of $C$ and let $\pi_i$ be a local equation of $C_i$. Let $K_i$ be the fraction field of the henselization of the local ring at the generic point $\xi_i$ of $C_i$.

 For a character $\chi\in\rmH^1(U,\Q/\Z)$, let $\chi|_{K_i}$ denote the pullback of $\chi$ by the map $\Spec K_i\to U$. We put $\Sw(\chi)=\sum_{i\in I} \sw(\chi|_{K_i}) C_i$, $\Dt(\chi)
 =\sum_{i\in I}\dt(\chi|_{K_i}) C_i$ and let $Z_\chi$ be the support of $\Sw(\chi)$. 
 We can uniquely define a global section
 \[
\rsw(\chi)\in \Gamma(Z_\chi, \Omega^1_X(\log C)(\Sw(\chi)))
\] 
whose stalk of $\rsw(\chi)$ at $\xi_i$ is equal to $\rsw(\chi|_{K_i})$ by \cite[Theorem 7.1, 7.3]{Ka89}. 
Similarly, we can can uniquely define a global section
\[
\ch'(\chi)\in \Gamma(Z_\chi, F\Omega^1_X(p\cdot\Dt(\chi)))
\] 
whose stalk of $\ch'(\chi)$ at $\xi_i$ is equal to $\ch'(\chi|_{K_i})$ by \cite[Theorem 1.5]{O25}. Here, $F\Omega^1_X$ denotes the Frobenius--Witt differential introduced by Saito \cite{Sa22}.

\begin{definition}
We say that $(X, U, \chi)$ is strongly clean at $x\in X$ if one of the following conditions is satisfied:
\begin{enumerate}[label=\textup{\arabic*.}]
\item $x\notin Z_{\chi}$. 

\item $x\in Z_{\chi}$ and if we write $\rsw \chi|_{D_i, x}=(\alpha_i+\beta_i d\log \pi_i)/\pi_i^{\sw(\chi|_{K_i})}$ where $\alpha_i\in \Omega^1_{\calO_{D_i, x}}, \beta_i\in \calO_{D_i,x}$, the image of $\beta_i$ in $k(x)$ is non-zero for every $i=1,\dots, r$. 
\end{enumerate}
We say that $(X, U, \chi)$ is strongly clean if it is strongly clean at every $x\in X$. 

\end{definition}
\begin{lemma}\label{1.6}
Assume that $\chi|_{K_i}$ is of type \textup{II} for every $i=1,\dots, r$. Then $(X,U, \chi)$ is not strongly clean. 
\end{lemma}
\begin{proof}
Under assumption, the refined Swan conductor $\rsw(\chi|_{K_i})$ is the image of the characteristic form $\ch(\chi|_{K_i})$ under the map $\rmH_1(L_{{F_i}^{1/p}/\calO_{K_i}})\to \Omega^1_{F_i}\otimes_{F_i} {F_i}^{1/p}\to \Omega^1_{F_i}(\log)\otimes_{F_i}{F_i}^{1/p}$ by \cite[Proposition 3.4]{O25}. Therefore, the coefficient of $d\log \pi_i$ in $\rsw(\chi|_{K_i})$ is zero for every $i=1,\dots , r$. 
\end{proof}

We prove some properties of the total dimension related to the blow-up and the specialization. These are non-logarithmic variants of the properties of the Swan conductor proved by Kato \cite{Ka89}. In the equal characteristic situation, these properties of the total dimension were proved by Matsuda \cite{Ma97}.

\begin{prop}\label{1.9}
Let $X=\Spec A$ be the spectrum of a regular excellent local ring $A$ with a perfect residue field. 
  Let $X'\to X$ be the blow-up at the closed point. Let $\nu$ be the generic point of the closed fiber. Then we have 
  \[
  \dt_{\nu}(\chi)\le \sum_{i=1}^r\dt_i(\chi). 
  \]
 Here, $\dt_{\nu}(\chi)$ denotes the total dimension of the pullback of $\chi$ to the local field at $\nu$ and $\dt_i(\chi)=\dt(\chi|_{K_i})$. 
\end{prop}
\begin{proof}
  By \cite[Theorem 8.1]{Ka89}, the inequality  
  \[
  \sw_{\nu}(\chi)\le \sum_{i=1}^r\sw_i(\chi)
  \]
  holds and the equality holds if and only if $\chi$ is strongly clean.
  If there exists $i$ such that $\chi|_{K_i}$ is of type I, we have 
  \[
  \dt_{\nu}(\chi)\le \sw_{\nu}(\chi)+1\le \sum_{i=1}^r\sw_i(\chi)+1\le \sum_{i=1}^r\dt_i(\chi).
  \]
  If $\chi|_{K_i}$ is of type II for every $i$, we have  
   \[
  \dt_{\nu}(\chi)\le  \sw_{\nu}(\chi)+1\le \sum_{i=1}^r\sw_i(\chi)= \sum_{i=1}^r\dt_i(\chi)
  \]
  because $\chi$ is not strongly clean in this case by Lemma \ref{1.6}. 
  \end{proof}
\begin{prop}\label{1.8}
 Let $X=\Spec A$ be the spectrum of a regular excellent local ring $A$ with a perfect residue field.  Let $\frakm$ be the maximal ideal of $A$. Assume that $r=1$ and let $\frakp$ be the prime ideal of $A$ defining $C_1$. Let $\frakq\in \Spec(A)_1$ and assume that $\frakm=\frakp+\frakq$. Let $\chi|_{\frakq}$ denote the pullback of $\chi$ by closed immersion $\Spec A/\frakq\to \Spec A$.  
 We put $n=\dt_{\frakp}(\chi), n'=\dt(\chi|_{\frakq})$ and 
 \[
 \ch_{\frakp}(\chi)=(\alpha w\pi_1+\beta)/\pi^n_1
  \]
  \[
 \ch(\chi|_{\frakq})=\alpha' w\pi_1/\pi^{n'}_1
  \]
 where $\alpha^p\in A/\frakp, \beta\in \Omega^1_{A/\frakp}$ and $\alpha'\in A/\frakm$.  
  
 Then, 
 
 \textup{(1)} We have 
  \[
  \dt(\chi|_{\frakq})\le \dt_{\frakp}(\chi).
  \]

  \textup{(2)} The equality holds if and only if one of the following conditions is satisfied.
  \begin{enumerate}[label=\textup{\arabic*.}]
  \item $(X,U,\chi)$ is strongly clean. 
  
  \item $\chi$ is of type \textup{II} at $\frakp$ and the image of $\alpha^p$ in $A/\frakm$ is non-zero.
  \end{enumerate}
  
  \textup{(3)} When the equality holds, the image of $\alpha^p$ in $A/\frakm$ is ${\alpha'}^p$.  
  \end{prop}

\begin{proof}
(1), (2) Since the residue field of $A/\frakq$ is perfect, the total dimension $\dt(\chi|_{\frakq})$ is equal to $\sw(\chi|_{\frakq})+1$. By \cite[Theorem 9.1]{Ka89}, we have 
  \[
  \sw(\chi|_{\frakq})\le \sw_{\frakp}(\chi)
  \]
  and the equality holds if and only if $(X, U,\chi)$ is strongly clean. 
  If $\chi$ is of type I at $\frakp$, we have 
  \[
  \dt(\chi|_{\frakq})= \sw(\chi|_{\frakq})+1\le \sw_{\frakp}(\chi)+1=\dt_{\frakp}(\chi). 
  \]
  The equality holds if $(X, U,\chi)$ is strongly clean. 
    
Next, consider the case where $\chi$ is of type II at $\frakp$. In this case, $(X, U,\chi)$ is not strongly clean by Lemma \ref{1.6}, so we have 
 \[
  \dt(\chi|_{\frakq})= \sw(\chi|_{\frakq})+1\le \sw_{\frakp}(\chi)=\dt_{\frakp}(\chi). 
  \]
  
We may assume $\dim A=2$ for the same reason as in the proof of \cite[Theorem 9.1]{Ka89}. Let $\tau_2$ be an element of $A$ satisfying $\frakq=(\tau_2)$ and let $A'=A[u_2]/(u_2^p-\tau_2)$. By \cite[Lemma 4.1.1]{O25}, using the notation of loc. cite., if the coefficient of $w\pi_1$ in $\ch(\chi|_{K_1})$ is non-zero, then we have $\dt(\chi'|_{K'_1})=\dt(\chi|_{K_1})$ and $\chi'|_{K'_1}$ is of type I. Hence, we have 
\[
\dt(\chi|_{K_1})=\dt(\chi'|_{K'_1})\ge \dt(\chi'|_{\frakq'})=\dt(\chi|_{\frakq})
\]
where $\frakq'=(u_2)\subset A'$. The equality holds if and only if $(X', U',\chi')$ is strongly clean. Since the coefficient of $d\log \pi_1$ in $\rsw(\chi'|_{K'_1})$ is $\alpha$, the triple $(X',U',\chi')$ is strongly clean if and only if the image of $\alpha$ in $(A/\frakm)^{1/p}$ is non-zero. If the coefficient of $w\pi_1$ in $\ch(\chi|_{K_1})$ is zero, we have 
\[
\dt(\chi|_{K_1})>\dt(\chi'|_{K'_1})\ge \dt(\chi'|_{\frakq'})=\dt(\chi|_{\frakq})
\]
and the equality does not hold. 

(3) If $\chi$ is of type I at $\frakp$, the characteristic form is the image of the refined Swan conductor and the refined Swan conductor $\rsw(\chi|_{\frakq})$ is the image of the $\rsw_{\frakp}(\chi)$ by \cite[Theorem 9.1]{Ka89}. Therefore, the claim follows. 

If $\chi$ is of type II at $\frakp$ and the equality holds, the coefficient of $d\log \pi_1$ in $\rsw(\chi'|_{K'_1})$ is $\alpha$ and applying \cite[Theorem 9.1]{Ka89}, the claim follows. 

\end{proof}

\begin{definition}
  Let $D=\sum_{i=1}^r m_i C_i$ be an effective Cartier divisor. We define $\fil_D \rmH^1(U, \Q/\Z)$ as the subgroup consisting of $\chi\in \rmH^1(U,\Q,\Z)$ satisfying $\dt(\chi|_{K_i})\le m_i$ for every $i=1,\dots, r$.   
  
  We define the subgroup $\pi^{\ab}_1(X, D)$ of the abelian etale fundamental group by
  \[
  \pi^{\ab}_1(X,D)=\Hom(\fil_D\rmH^1(U,\Q/\Z), \Q/\Z). 
  \]
   \end{definition}
   
    The definition of the abelian fundamental group with modulus of a smooth scheme over a perfect field is given in \cite{KS14}.
    
  \begin{definition}
We define $\fil_D^{\log} \rmH^1(U, \Q/\Z)$ as the subgroup consisting of $\chi\in \rmH^1(U,\Q,\Z)$ satisfying $\sw(\chi|_{K_i})\le m_i$ for every $i=1,\dots, r$. 
\end{definition}  

We give an alternative definition of the abelian fundamental group with modulus when $X$ is a regular flat scheme of finite type over a discrete valuation ring of mixed characteristic $(0,p)$ with perfect residue field. 

\begin{thm}\label{1.2}
 Under the assumption on $X$, let $D=\sum_{i=1}^r m_i C_i$ be an effective Cartier divisor. For $\chi\in \rmH^1(U, \Q/\Z)$, the following conditions are equivalent:
  
  \textup{(1)} For every $i=1, \dots, r$, we have $\dt(\chi|_{K_i})\le m_i$. 

  \textup{(2)} For an integral flat curve $Z\subset X$ not contained in $C$, let $\phi_{Z}:Z^N\to Z$ be the normalization map. For all such curves $Z$, the character $\chi|_Z\in \rmH^1(\phi_Z^{-1}(U), \Q/\Z)$ satisfies 
  \begin{equation}\label{1.10}
  \sum_{y\in Z_{\infty}} \dt_y(\chi|_Z)\cdot [y]\le \phi^*_{Z} D
  \end{equation}
  where $Z_{\infty}\subset Z^N$ is the finite set of the points $y$ such that $\phi_Z(y)\notin U$. The total dimension $\dt_y(\chi|_Z)$ means $\sw_y(\chi|_Z)+1$. 
\end{thm}
\begin{proof}
We prove the implication $(1)\Rightarrow(2)$. For an integral curve $Z$, there exists a chain of blow-ups at closed points $f: X'\to X$ such that the strict transform $Z'$ of $Z$ is regular and $Z'$ intersects with $C$ transversally by \cite[Theorem A.1]{Ja10}. Let $\chi'$ be the pullback of $\chi$ to $U'=U\times_X X'$. By Proposition \ref{1.9}, we have $\chi'\in \fil_{f^*D}\rmH^1(U, \Q/\Z)$. By Proposition \ref{1.8} (1), we have 
\[
\sum_{y\in Z'_{\infty}}\dt_y(\chi|_{Z'})\cdot [y]\le \phi^*_{Z'}f^*D. 
\]
This inequality implies the inequality \eqref{1.10}. 

We prove the implication $(2)\Rightarrow (1)$. We may assume $r=1$. The set of closed points on which neither of the conditions in Proposition \ref{1.8} (2) is satisfied is a closed subset of codimension $\ge 2$ in $X$. Therefore, there exist a curve $Z$ and a point $y\in Z_{\infty}$ such that we have $\dt(\chi|_{K_1})=\dt_y(\chi|_Z)$. 
\end{proof}

\begin{lemma}\label{1.7}
Let $\map{f}{Y}{X}$ be a map of regular excellent schemes. Assume that $C'=f^{-1}(C)$ is a reduced simple normal crossing divisor on $Y$ and that $C'_i=f^{-1}(C_i)$ is regular irreducible for every $i=1,\dots,r$. Let $E$ be the pullback of $D$ to $Y$ and $V$ be the complement of $E$. Then there exist the following commutative diagrams:
\[
 \xymatrix{
\fil_D^{\log}\rmH^1(U,\Q/\Z)/\fil_{D-C_i}^{\log}\rmH^1(U,\Q/\Z)\ar[r]^-{\rsw}\ar[d]& \rmH^0(C_i, \Omega^1_X(\log C)(D)|_{C_i})\ar[d] \\
\fil_E^{\log}\rmH^1(V,\Q/\Z)/\fil_{E-C'_i}^{\log}\rmH^1(V, \Q/\Z)\ar[r]^-{\rsw}& \rmH^0(C'_i, \Omega^1_Y(\log C')(E)|_{C'_i})
}  
\]

\[
 \xymatrix{
\fil_D\rmH^1(U,\Q/\Z)/\fil_{D-C_i}\rmH^1(U,\Q/\Z)\ar[r]^-{\ch'}\ar[d]& \rmH^0(C_i, F\Omega^1_X(D)|_{C_i})\ar[d] \\
\fil_E\rmH^1(V,\Q/\Z)/\fil_{E-C'_i}\rmH^1(V, \Q/\Z)\ar[r]^-{\ch'}& \rmH^0(C'_i, F\Omega^1_Y(E)|_{C'_i})
}  
\]
\end{lemma}
\begin{proof}
It suffices to prove in the cases where $f$ is smooth and $f$ is a closed immersion. 
If $f$ is smooth, the assertion follows from the functoriality and integrality of the refined Swan conductor and the characteristic form. If $f$ is a closed immersion, the assertion follows from \cite[Theorem 9.1]{Ka89} and Proposition \ref{1.8} (3). 

\end{proof}

\section{Main Theorem}
Let $K$ be a complete discrete valuation field of characteristic zero with perfect residue field $k$ of characteristic $p>0$. Let $X$ be a connected strictly semi-stable projective scheme over $\calO_K$. Let $X_K$ and $X_k$ be the generic fiber and the closed fiber, respectively. 

Let $C$ be a reduced simple normal crossing divisor containing the closed fiber $X_k$ and let $\{C_j\}_{j\in I}$ be the irreducible components. 
Let $D=\sum_{j\in I} m_j C_j$ be an effective Cartier divisor with support $C$. 
Let $K_j$ be the fraction field of the completion of the local ring of $X$ at the generic point of $C_j$.  Let $e_j$ be the ramification index of $K_j$ and write $e'_j=\lceil e_jp/(p-1)\rceil$. 

Let $Y$ be a regular hypersurface section flat over $\calO_K$ such that $Y\times_X C$ is a reduced simple normal crossing divisor and $Y\times_X C_j$ is regular irreducible for any $j\in I$. Let $\map{i}{Y}{X}$ be the closed immersion. Let $E$ be the simple normal crossing divisor $Y\times_X D$. 
We fix some notation. Let $D_1, D_e, D_{e'}$ be the divisors $\sum_{j\in I}C_j, \sum_{j\in I}e_jC_j, \sum_{j\in I}e'_jC_j$,  respectively. Let $E_1, E_e, E_{e'}$ be the pullbacks to $Y$, respectively. 

We note that if we let $K'_j$ be the fraction field of the completion of the local ring of $Y$ at the generic point of $C_j\cap Y$, the absolute ramification index of $K_j$ and $K'_j$ is equal. 

Let $\tilde{I}$ be the subset of $I$ consisting of $j$ such that $m_j< e'_j+e_j$ and $C_j$ is contained in the closed fiber $X_k$. We put
\[
\tilde{D}=\sum_{j\in \tilde{I}}(e'_j+e_j)C_j+\sum_{j\notin \tilde{I}}m_jC_j.
\]

Let $I'$ be the subset of $I$ consisting of $j$ such that $p$ divides $m_j$ and $C_j$ is contained in the closed fiber $X_k$. We put
\[
D'=\sum_{j\in I'}(m_j+1)C_j+\sum_{j\notin I'}m_jC_j. 
\]
\begin{definition}\label{2.2}
We say that $Y$ is sufficiently ample for $(X,D)$ if the following conditions hold:
\begin{enumerate}[label=\textup{\arabic*.}]
  \item For any $j\in I$ such that $C_j\subset X_k$ and for any effective Cartier divisor $\Theta\le D_{e'}$, we have 
  \begin{equation}\label{2.6}
  \begin{aligned}
  &\rmH^0(C_j, \Omega^1_X(\log C)(\Theta-Y)|_{C_j})=0\\
 &\rmH^0(C_j, \calO_{C_j}(\Theta-Y)|_{C_j})=0\\
 & \rmH^1(C_j, \calO_{C_j}(\Theta-2Y)|_{C_j})=0  
    \end{aligned}
    \end{equation}
    
    \item For any $j\in \tilde{I}$ and for any effective Cartier divisor $\Theta\le \tilde{D}$, we have three equalities \eqref{2.6}. 
    
  \item  For any $j\in I'$, we have
  \[
  \rmH^0(C_j, F\Omega^1_X(D'-Y)|_{C_j})=0\\
  \]
    \end{enumerate}

\end{definition}
\begin{remark}
  When $\dim X_k\ge 2$, there exists an integer $N$ such that if the degree of $Y$ is greater than $N$, then $Y$ is sufficiently ample for $(X,D)$. 
\end{remark}
\begin{thm}\label{main}
Assume that $Y$ is sufficiently ample for $(X,D)$. Then, the natural map 
\[
\pi_1^{\ab}(Y,E)\lra \pi_1^{\ab}(X,D)
\]
is an isomorphism when $\dim X_K\ge 3$ and a surjection when $\dim X_K=2$.

\begin{lemma}\label{2.1}
Let $U$ be the complement of $D$.  The restriction map
  \[
i^*:\rmH^1(U,\Q/\Z) \lra \rmH^1(U\cap Y, \Q/\Z)
\]
is an isomorphism when $\dim X_K\ge 3$ and is an injection when $\dim X_K=2$.  
\end{lemma}
\begin{proof}
We apply \cite[Theorem 1.1]{EK16}. Since $X_K$ is connected projective over $K$ and $Y_K$ is a regular hypersurface section transversal to $D_K$, the condition $\mathrm{Leff}(X_K, Y_K)$ is satisfied when $\dim X_k\ge 3$ and $\mathrm{Lef}(X_K, Y_K)$ is satisfied when $\dim X_K=2$ by \cite[Expose X. Exemple 2.2.]{SGA2}. 
Since the characteristic of $K$ is zero, the natural map 
\[
\pi_1(Y_K\backslash E_K)\lra \pi_1(X_K\backslash D_K)
\]
is an isomorphism when $\dim X_K=3$ and a surjection when $\dim X_K=2$. 
As we assume that $D$ contains the closed fiber $X_k$, we have $Y_K\backslash E_K=U\cap Y$ and $X_K\backslash D_K=U$. Hence, the assertion follows. 
\end{proof}

\end{thm}
It suffices to show Theorem \ref{tame} and Theorem \ref{nonlog}. For an abelian group $M$, let $M\{p'\}$ be the prime-to-$p$ torsion part of $M$. 
\begin{thm}\label{tame}
The natural map 
\[
\fil_D \rmH^1(U,\Q/\Z)\{p'\}\lra \fil_E\rmH^1(U\cap Y,\Q/\Z)\{p'\}
\]
is an isomorphism when $\dim X_K\ge 3$ and an injection when $\dim X_K=2$.
\end{thm}
\begin{proof}
For an integer $n$ prime to $p$, 
\[
\fil_D \rmH^1(U,\Z/n\Z)=\rmH^1(U,\Z/n\Z)
, \ \fil_E \rmH^1(U\cap Y,\Z/n\Z)=\rmH^1(U\cap Y,\Z/n\Z)
\]
Therefore, the isomorphism follows from \cite[Theorem 1.1]{EK16}. 
\end{proof}
Before we prove Theorem \ref{nonlog}, we prove the logarithmic case. 

\begin{thm}\label{log}
Assume that $Y$ is sufficiently ample for $(X,D)$. Then, the natural map 
\[
\fil_D^{\log} \rmH^1(U,\Z/p^n\Z)\lra \fil_E^{\log}\rmH^1(U\cap Y,\Z/p^n\Z)
\]
is an isomorphism for all $n\ge 1$ when $\dim X_K\ge 3$ and an injection when $\dim X_K=2$.
\end{thm}

By Lemma \ref{2.1}, it suffices to show that the map is surjective. 
We first prove the case $n=1$. 
\begin{proof}[Proof $(n=1)$]
If we have $D\ge D_{e'}$, we have equalities
\[
\fil_D^{\log}\rmH^1(U,\Z/p\Z)=\rmH^1(U,\Z/p\Z), 
\]
\[
\fil_E^{\log}\rmH^1(U\cap Y, \Z/p\Z)=\rmH^1(U\cap Y,\Z/p\Z)
\]
by Proposition \ref{1.3}. Hence, the map
\[
\fil_D^{\log}\rmH^1(U,\Z/p\Z)\lra \fil_E^{\log}\rmH^1(U\cap Y,\Z/p\Z)
\]
is an isomorphism by Lemma \ref{2.1}. 

It suffices to show that if the theorem holds for $\Theta=\sum_{j\in I}n_j C_j$, then the theorem holds for $\Theta-C_j$ when $n_j\ge 2$. 
We write $\Xi=\Theta\times_X Y$. 

We note that if the divisor $C_j$ intersects the generic fiber $X_K$, for every character $\chi\in H^1(U, \Q/\Z)$, the restriction $\chi|_{K_j}$ is tame. Hence, it suffices to consider the case  $C_j\subset X_k$. 

We have to show that the map 
\[
\fil_{\Theta-C_j}^{\log}\rmH^1(U,\Q/\Z)\lra \fil_{\Xi-C_j\cap Y}^{\log}\rmH^1(U\cap Y,\Q/\Z)
\]
is surjective. By the surjectivity of 
\[
\fil_{\Theta}^{\log}\rmH^1(U,\Z/p\Z)\lra \fil_{\Xi}^{\log}\rmH^1(U\cap Y,\Z/p\Z),
\]
it suffices to show that the map 
\[
\fil_{\Theta}^{\log}\rmH^1(U,\Z/p\Z)/\fil_{\Theta-C_j}^{\log}\rmH^1(U,\Z/p\Z)\lra \fil_{\Xi}^{\log}\rmH^1(U\cap Y,\Z/p\Z)/\fil_{\Xi-C_j\cap Y}\rmH^1(U\cap Y, \Z/p\Z)
\]
is injective. We have the following commutative diagram
  \begin{equation}\label{2.4}
  \vcenter{\hbox{
  \scalebox{0.9}{$
  \xymatrix{
\fil_{\Theta}^{\log}\rmH^1(U,\Z/p\Z)/\fil_{\Theta-C_j}^{\log}\rmH^1(U,\Z/p\Z)\ar[r]^-{\rsw}\ar[d]& \rmH^0(C_j, \Omega^1_X(\log C)(\Theta)|_{C_j})\ar[d] \\
\fil_{\Xi}^{\log}\rmH^1(U\cap Y,\Z/p\Z)/\fil_{\Xi-C_j\cap Y}\rmH^1(U\cap Y, \Z/p\Z)\ar[r]^-{\rsw}& \rmH^0(C_j\cap Y, \Omega^1_Y(\log C\cap Y)(\Xi)|_{C_j})
} $ 
}
}
}
  \end{equation}
   
   Since the refined Swan conductor is injective, it suffices to show that the right vertical map is injective. 
  There exists the following commutative diagram of exact sequences:
  \[
  \xymatrix{
  0\ar[r]& \Omega^1_X(\log C)(-Y)\ar[r]\ar[d] &\Omega^1_X(\log C) \ar[r]\ar[d]& \Omega^1_X(\log C)\otimes_{\calO_{X}}\calO_{Y}\ar[r]\ar[d]&0\\
  0\ar[r] &0\ar[r] &i_*\Omega^1_Y(\log C\cap Y)\ar[r] &i_*\Omega^1_Y(\log C\cap Y)\ar[r]&0
  }\]
  If $\mathcal{K}$ denotes the kernel of the middle map, we get an exact sequence
  \[
  0\lra \Omega^1_X(\log C)(-Y)\lra \mathcal{K}\lra \calO_X(-Y)\otimes \calO_Y\lra 0. 
  \]
  Hence, we have an exact sequence
   \[
  0\lra \Omega^1_X(\log C)(\Theta-Y)\otimes \calO_{C_j}\lra \mathcal{K}(\Theta)\otimes \calO_{C_j}\lra \calO_X(\Theta-Y)\otimes \calO_{C_j\cap Y}\lra 0. 
  \]
  
  We also have an exact sequence 
 \[
 0\lra \calO_{C_j}(\Theta-2Y)\lra \calO_{C_j}(\Theta-Y)\lra \calO_{C_j}(\Theta-Y)\otimes\calO_{Y\cap C_j}\lra 0.
 \]
 From the condition Definition \ref{2.2}.1, we have 
 \[
 \rmH^0(C_j, \mathcal{K}|_{C_j})=0. 
 \]
 Therefore, the right vertical map in the diagram \eqref{2.4} is injective and the assertion follows.

\end{proof}

Next, we prove the case $n\ge 2$. 

\begin{proof}[Proof $(n\ge 2)$]
  We prove the claim by induction on $n$. By Lemma \ref{2.1}, there exists the inverse map 
  \[
(i^*)^{-1}:\rmH^1(U\cap Y, \Z/p^m\Z)\lra \rmH^1(U,\Z/p^m\Z)
  \]
 for any positive integer $m$.
 
 Assume $D\ge D_{e'}+D_e$. Take $\chi \in \fil_E^{\log}\rmH^1(U\cap Y, \Z/p^n\Z)$ such that $\Sw(\chi)> E_{e'}$. Then, by Lemma \ref{1.4}, we have $\Sw(p\chi)=\Sw(\chi)-E_e$. By the induction hypothesis, we have $\Sw((i^*)^{-1}(p\chi))\le D-D_e$. 
Since $i^*$ and $p\cdot $ is commutative, again by Lemma \ref{1.4}, we have $\Sw((i^*)^{-1}(\chi))\le D$.

Take $\chi \in \fil^{\log}_E\rmH^1(U\cap Y, \Z/p^n\Z)$ such that we do not have $\Sw(\chi)> E_{e'}$. Put $\psi=(i^*)^{-1}\chi$. Assume $\Sw(\psi)>D$. By Lemma \ref{1.4}, we have $\Sw(p\psi)>D-D_e$. By the induction hypothesis, we have $\Sw(pi^*\psi)=\Sw(p\chi)>E-E_e$. By Lemma \ref{1.4}, we have $\Sw(\chi)>E-E_{e}\ge E_{e'}$ and this is a contradiction.

Hence, under the assumption $D\ge  D_{e'}+D_e$, the map $\fil_D^{\log}\rmH^1(U,\Z/p^n\Z)\to \fil_E^{\log}\rmH^1(U\cap Y,\Z/p^n\Z)$ is surjective.

In the case where we do not have $D\ge D_{e'}+D_e$, the proof is similar to the case $n=1$. Let $\Theta=\sum_{j\in I} n_j C_j\le \tilde{D}$ be an effective Cartier divisor with support $C$ and $\Xi=\Theta\times_X Y$. For $j\in \tilde{I}$ such that $n_j\ge 2$, consider the following commutative diagram
\begin{equation*}
    \xymatrix{
\fil_{\Theta}^{\log}\rmH^1(U,\Z/p^n\Z)/\fil_{\Theta-C_j}^{\log}\rmH^1(U,\Z/p^n\Z)\ar[r]^-{\rsw}\ar[d]& \rmH^0(C_j, \Omega^1_X(\log C)(\Theta)|_{C_j})\ar[d] \\
\fil_{\Xi}^{\log}\rmH^1(U\cap Y,\Z/p^n\Z)/\fil^{\log}_{\Xi-C_j\cap Y}\rmH^1(U\cap Y, \Z/p^n\Z)\ar[r]^-{\rsw}& \rmH^0(C_j\cap Y, \Omega^1_Y(\log C\cap Y)(\Xi)|_{C_j\cap Y})
}
  \end{equation*}
Using the equalities in Definition \ref{2.2}.2, the right vertical map is injective by the same argument as in the case $n=1$, so the left vertical map is also injective. Therefore, if the map
\[
\fil_{\Theta}^{\log}\rmH^1(U,\Z/p^n\Z)\lra \fil_{\Xi}^{\log}\rmH^1(U\cap Y,\Z/p^n\Z)
\]
is surjective, then the map
\[
\fil_{\Theta-C_j}^{\log}\rmH^1(U,\Z/p^n\Z)\lra \fil_{\Xi-C_j \cap Y}^{\log}\rmH^1(U\cap Y,\Z/p^n\Z)
\]
is surjective. The assertion follows by descending induction on $n_j$ for each $j\in \tilde{I}$, starting from $\tilde{D}$.

\end{proof}

We deduce the non-logarithmic case from the logarithmic case. The proof essentially follows the same steps as in \cite{KS14}. 

\begin{thm}\label{nonlog}
Assume that $Y$ is sufficiently ample for $(X,D)$. Then, the natural map 
\[
\fil_D \rmH^1(U,\Z/p^n\Z)\lra \fil_E\rmH^1(U\cap Y,\Z/p^n\Z)
\]
is an isomorphism for all $n\ge 1$ when $\dim X_K\ge 3$ and an injection when $\dim X_K=2$.
\end{thm}
\begin{proof}
The injectivity follows from Lemma \ref{2.1}. Thus, it suffices to show the surjectivity when dim $X_K\ge 3$. 
Noting that each coefficient of $C_j\subset X_k$ in $D'$ is not divided by $p$, by Corollary \ref{1.11}, we have 
\begin{align*}
\fil_{D'}H^1(U,\Z/p^n\Z)&=\fil^{\log}_{D'-C}H^1(U, \Z/p^n \Z)\\
\fil_{E'}H^1(U\cap Y, \Z/p^n \Z)&=\fil^{\log}_{E'-C\cap Y}H^1(U\cap Y, \Z/p^n\Z). 
\end{align*}
By Theorem \ref{log}, the map 
\[
\fil_{D'} \rmH^1(U,\Z/p^n\Z)\lra \fil_{E'}\rmH^1(U\cap Y,\Z/p^n\Z)
\]
is an isomorphism. 
So, it suffices to show that 
\[
\fil_{D'}\rmH^1(U,\Z/p^n\Z)/\fil_{D}\rmH^1(U, \Z,p^n\Z)\lra \fil_{E'}\rmH^1(U\cap Y,\Z/p^n\Z)/\fil_E\rmH^1(U\cap Y, \Z/p^n\Z)
\]
is injective. 
Consider the following commutative diagram:
\begin{equation}
\vcenter{\hbox{
  \scalebox{0.9}{$
    \xymatrix{
\fil_{D'}\rmH^1(U,\Z/p^n\Z)/\fil_{D}\rmH^1(U, \Z/p^n\Z)\ar[r]^-{\ch'}\ar[d]& \bigoplus_{j\in I'}\rmH^0(C_j, F\Omega^1_X(D')|_{C_j})\ar[d] \\
\fil_{E'}\rmH^1(U\cap Y,\Z/p^n\Z)/\fil_E\rmH^1(U\cap Y, \Z/p^n\Z)\ar[r]^-{\ch'}& \bigoplus_{j\in I'}\rmH^0(C_j\cap Y, F\Omega^1_Y(E')|_{C_j\cap Y})
}
$}
}
}
  \end{equation}
  Since the characteristic form is injective, it suffices to show that the right vertical map is injective. 
  There exists following commutative diagram of exact sequences:
  \begin{equation}\label{2.3}
  \xymatrix{
  0\ar[r]& F\Omega^1_X(-Y_k)\ar[r]\ar[d] &F\Omega^1_X \ar[r]\ar[d]& F\Omega^1_X\otimes_{\calO_{X_k}}\calO_{Y_k}\ar[r]\ar[d]&0\\
  0\ar[r] &0\ar[r] &i_*F\Omega^1_Y\ar[r] &i_*F\Omega^1_Y \ar[r]&0
  }
  \end{equation}
  By \cite[Corollary 2.8]{Sa22}, the kernel of the right vertical map is equal to $\rmF^*(\calO_X(-Y)\otimes_{\calO_X} \calO_{Y_k})$ where $\rmF: Y_k\to Y_k$ is the Frobenius map. 
  Let $\mathcal{L}$ be the kernel of the middle map. Then we get an exact sequence
    \[
  0\lra F\Omega^1_X(-Y_k)\lra \mathcal{L} \lra \rmF^*(\calO_X(-Y)\otimes_{\calO_X} \calO_{Y_k})\lra 0
  \]
  and we get an exact sequence
   \[
  0\lra F\Omega^1_X(D'-Y_k)\otimes \calO_{C_j}\lra \mathcal{L}(D')\otimes\calO_{C_j}\lra \rmF^*(\calO_{C_j}(D'-Y)\otimes_{\calO_{C_j}} \calO_{Y_k\cap C_j})\lra 0. 
  \]
  
 We also have an exact sequence 
 \[
 0\lra \calO_{C_j}(D'-2Y)\lra \calO_{C_j}(D'-Y)\lra \calO_{C_j}(D'-Y)\otimes\calO_{Y_k\cap C_j}\lra 0.
 \]
 From the equality in Definition \ref{2.2}.3, we have 
 \[
 \rmH^0(C_j, \mathcal{L}|_{C_j})=0. 
 \]
 Therefore, the right vertical map in the diagram \eqref{2.3} is injective and the assertion follows. 
 
  \end{proof}
  
  \begin{cor}
Let X be a regular strictly semi-stable projective scheme over $\calO_K$. Assume that $Y$ is sufficiently ample for $(X,D)$. Then, the natural map 
\[
\fil_D \rmH^1(U,\Q/\Z)\lra \fil_E\rmH^1(U\cap Y,\Q/\Z)
\]
is an isomorphism when $\dim X_K\ge 3$ and an injection when $\dim X_K=2$.

Dually, the natural map
\[
\pi_1^{\ab}(Y, E)\lra \pi_1^{\ab}(X, D)
\]
is an isomorphism when $\dim X_K\ge 3$ and a surjection when $\dim X_K=2$.

\end{cor}

\textsc{Graduate School of Mathematical Sciences, University of Tokyo, 3-8-1 Komaba, Meguro-Ku, Tokyo 153-8914, Japan}

\textit{Email address}: \texttt{ooe-ryosuke075@g.ecc.u-tokyo.ac.jp}

\end{document}